\documentclass[12pt, twoside, a4paper]{article}
\usepackage[T1]{fontenc}
\usepackage{makeidx}\usepackage{amssymb}
\usepackage{amsmath,amssymb,amsopn,amsfonts,mathrsfs,amsbsy,amscd}
\usepackage[usenames,dvipsnames]{pstricks} \usepackage{epsfig}

\newcommand{\prs}{\langle\;,\;\rangle}
\newcommand{\too}{\longrightarrow}
\newcommand{\om}{\omega}
\newcommand{\esp}{\quad\mbox{and}\quad}

\newcommand{\G}{{\mathfrak{g}}}

\newcommand{\h}{{\mathfrak{h}}}
\newcommand{\ad}{{\mathrm{ad}}}

\newcommand{\tr}{{\mathrm{tr}}}

\newcommand{\B}{{\cal B}}

\newcommand{\Om}{\Omega}
\newcommand{\na}{\nabla}

\newcommand{\al}{\alpha}
\newcommand{\be}{\beta}
\newcommand{\ga}{\gamma}
\newcommand{\Ga}{\Gamma}

\newcommand{\de}{\delta}

\font\bb=msbm10

\def\K{\hbox{\bb K}}

\def\B{\hbox{\bb B}}
\def\R{\hbox{\bb R}}

\def\N{\hbox{\bb N}}

\def\C{\hbox{\bb C}}
\newtheorem{Def}{Definition}[section]
\newtheorem{theo}{Theorem}[section]
\newtheorem{pr}{Proposition}[section]
\newtheorem{Le}{Lemma}[section]
\newtheorem{co}{Corollary}[section]

\newtheorem{exem}{Example}
\newtheorem{rem}{Remark}

\title{ Special bi-invariant linear connections on Lie groups and finite dimensional
Poisson
structures
}
\author{Sa\"id Benayadi and Mohamed Boucetta
\footnote{This research was conducted within the framework of Action concert\'ee CNRST-CNRS Project SPM04/13.}
}\date{$\;$}

\begin{document}\maketitle

\begin{abstract}Let $G$ be a connected Lie group and $\G$ its Lie algebra. We denote by
$\na^0$ the torsion free
bi-invariant linear connection on $G$  given by
$\na^0_XY=\frac12[X,Y],$ for any left invariant vector fields $X,Y$. A Poisson structure
on $\G$ is a commutative and associative product on $\G$ for which
$\ad_u$ is a derivation, for any $u\in\G$.
 A torsion free bi-invariant linear connections on $G$ which have the same
curvature as $\na^0$ is called special.
We show that there is a
bijection between the space of special connections on $G$ and the space of Poisson
structures
on  $\G$. We compute the holonomy Lie algebra of a special connection and we show that
the Poisson structures associated to special connections which have the same holonomy Lie
algebra as $\na^0$ possess interesting properties.
Finally, we study Poisson structures  on a Lie algebra and
we give a large class of examples which gives, of course, a large class of special
connections.
\end{abstract}

{\it Key words:}  Lie groups, Lie algebras, bi-invariant linear connections, Poisson algebras, symplectic Lie algebras, symplectic double extension, semi-symmetric linear connections.

{\bf MSC.}  17A32, 17B05, 17B30, 17B63, 17D25, 53C05.
\section{Introduction}\label{section1}
All vector spaces, algebras, etc. in this paper will be over a ground
field $\K$ of characteristic 0.\\
A {\it Poisson algebra} is a finite dimensional Lie algebra $(\G,[\;,\;])$
endowed
with a commutative and associative product $\circ$ such that, for any $u,v,w\in\G$,
\begin{equation}
 \label{eq0}[u,v\circ w]=[u,v]\circ w+v\circ [u,w].
\end{equation}
An algebra $(\mathrm{A},.)$ is called {\it Poisson admissible } if
$(\mathrm{A},[\;,\;],\circ)$ is a  Poisson algebra, where
\begin{equation}\label{oh}
[u,v]=u.v-v.u\esp u\circ v=\frac12(u.v+v.u).\end{equation}
Poisson algebras constitute an interesting topic in algebra and were studied by many
authors (see for instance \cite{GozRem,  MarRem, Shestakov}). This paper aims
to give some new insights on them  based on an interesting geometric
interpretation of these structures when the field is either $\R$ or $\C$ (see Theorem
\ref{main}). Let us present
briefly this geometric interpretation. \\ Let
  $G$  be  a Lie group with $\G=T_eG$ its Lie algebra. The linear connection $\na^0$
given by $\na^0_XY=\frac12[X,Y]$, where $X,Y$ are left invariant vector fields, is torsion
free, bi-invariant, complete and its curvature $K^0$ is given by
$K^0(X,Y)=-\frac14\ad_{[X,Y]}$. Moreover, $\na^0K^0=0$ and the holonomy Lie algebra of
$\na^0$ at $e$ is  $\mathfrak{h}^0=\ad_{[\G,\G]}$. The main fact (see  Section
\ref{section2}) is that there is a bijection
between the set of Poisson structures on  $\G$ and the space of
bi-invariant torsion free linear connections  on $G$ which have the same curvature
as $\na^0$. We call such
connections {\it special}. Moreover, we show  that any special connection is semi-symmetric, i.e., its  curvature tensor $K$ satisfies $K.K=0$ (see Proposition \ref{mainpr}). In general, the
holonomy Lie
algebra $\mathfrak{h}$ of a bi-invariant linear connection  is difficult to
compute,
however, we show that,
for a special connection, $\mathfrak{h}$  contains
$\mathfrak{h}^0$  and can be easily computed (see Lemma \ref{lemma2}). A special
connection whose holonomy Lie algebra coincides with
$\mathfrak{h}^0$ will be called {\it strongly special}. So, according to the bijection
above, to any real Poisson algebra
corresponds a unique special connection on any associated Lie group. Poisson algebras
whose corresponding special connection is strongly special are particularly interesting.
We call such Poisson algebras {\it strong}.
With this interpretation
in mind, we devote
 Section \ref{sectional} to the
study of the general properties of Poisson algebras and Poisson admissible algebras and we
give
some general methods to build new Poisson algebras from old ones (see Theorem
\ref{theostrange}).
 We show that
any symmetric Leibniz algebra is a strong Poisson admissible algebra
and the curvature of the corresponding special connection is parallel (See Theorem
\ref{theoleibniz}). By using  the  geometric interpretation  of Poisson
structures, we get a large class of Lie groups which carry a bi-invariant connection
 $\na$ (different from $\na^0$) which has the same curvature and the same linear holonomy
as $\na^0$ and
moreover the curvature of $\na$ is parallel. We get hence interesting examples of
connections with parallel torsion and curvature. Such connections were studied by
Nomizu \cite{nomi}.
 Recall that
symmetric Leibniz algebras constitute a subclass of Leibniz algebras  introduced by Loday
in
\cite{loday}.    At the end of Section \ref{sectional},
we show that there is no non trivial Poisson structure on  a semi-simple Lie algebra (see
Theorem \ref{semisimple}). This result
 generalizes a result by \cite{GozRem}. In Section \ref{ass}, we show that an
associative algebra is Poisson admissible if and only if the underline Lie algebra is
2-nilpotent and
 we give a description of associative Poisson admissible algebras which permit to build
many
examples. Section \ref{sectionsp} is devoted to the study of symplectic Poisson
algebras. It is well-known  that if $(\G,\om)$ is a symplectic Lie
algebra there is a product $\al^{\mathrm{a}}$ on $\G$ which is Lie-admissible and left
symmetric.
When the Lie algebra is real, $\al^{\mathrm{a}}$ defines a  left invariant flat torsion
free linear
connection $\na^{\mathrm{a}}$ on any associated Lie group $G$. By using the general method
to build
a torsion free  symplectic connection from any  torsion free connection introduced in
\cite{cahen}, we get from $\na^{\mathrm{a}}$  a  left invariant torsion
free connection $\na^{\mathrm{s}}$ for which  the left invariant symplectic form
associated to $\om$ is
parallel. To our knowledge this connection has never been considered before.
From $\na^{\mathrm{s}}$ we get a product $\al^{\mathrm{s}}$ on $\G$. We
show that $(\G,\al^{\mathrm{a}})$ is Poisson admissible  iff
$(\G,\al^{\mathrm{s}})$ is
Poisson admissible  and
this is equivalent to $\G$ is 2-nilpotent Lie algebra and $[\ad_u,\ad_v^*]=0$ for any
$u,v\in\G$ where $\ad_u^*$ is the adjoint of $\ad_u$ with respect to $\om$.
A symplectic
Lie algebra satisfying these conditions is called {\it symplectic Poisson algebra}. We
show that the symplectic double extension process introduced in \cite{DarMed} permits the
construction of all
symplectic Poisson algebras.  Lie groups whose Lie algebras are symplectic Poisson
possesses an important geometric property (see Theorem \ref{polynomial} and the following
remarks).
In Section \ref{metri}, we study the problem of metrizability  of special connections.
Indeed, we consider a real Lie algebra
$(\G,[\;,\;],\prs)$ endowed with a nondegenerate symmetric bilinear metric. We denote by
$\ell$ the Levi-Civita product associated to  $(\G,[\;,\;],\prs)$. We show that if $\prs$
is positive definite $(\G,\ell)$ is  Poisson admissible  iff $\prs$ is
bi-invariant and in
this case the associated Poisson product $\circ$ is trivial. We give a description of
$(\G,[\;,\;],\prs)$ for which $(\G,\ell)$ is  Poisson admissible  in the case
where
$[\G,\G]$ is nondegenerate and $\prs$ has any signature.

\section{Geometric interpretation of finite dimensional Poisson
structures}\label{section2}
We give in this section an interesting geometric
interpretation of Poisson structures involving the theory of connections and holonomy
algebras. This theory
 is a fundamental topic in differential
geometry and has its origin in the work of Elie Cartan \cite{cartan1, cartan4}. For a
detailed account of this
theory, see Ehresmann \cite{ehres},
Chern \cite{chern}, Lichnerowicz \cite{lichner}, Nomizu \cite{nomi}, and Kobayashi
\cite{ko}. Let us recall some classical facts about linear connections and state
some formulas which will lead naturally to the desired interpretation.\\
 Given a
linear connection on a
smooth manifold $M$, we consider the covariant differentiation $\na$ associated to it.
Let $T^\na$ and $K^\na$ be,
respectively, the torsion and curvature tensor fields on $M$ with respect to $\na$:
$$T^\na(X,Y)=\na_XY-\na_YX-[X,Y]\esp K^\na(X,Y)=
[\na_X,\na_Y]-\na_{[X,Y]}.$$
For any closed curve $\tau$ at $p\in M$, the parallel displacement along
$\tau$ is a linear transformation of $T_pM$, and the totality of these
linear transformations for all closed curves forms the {\it holonomy
group} $H(p)$. The {\it restricted holonomy group} $H(p)^0$ is the subgroup consisting of
parallel displacements along all closed curves which are homotopic to
zero. Its Lie algebra is called {\it holonomy Lie algebra}.
On the other hand,   consider linear
endomorphisms
of $T_pM$ of the form $K^\na(X, Y)$, $(\na_ZK^\na) (X, Y)$, $(\na_W\na_ZK^\na)
(X, Y)$, . . . (all covariant derivatives), where $X, Y, Z, W,$ . . . are arbitrary
tangent vectors at $p$. All these linear endomorphisms span a subalgebra
$\mathfrak{h}_p^\na$ of the Lie algebra consisting of all linear endomorphisms of $T_pM$.
We call it the {\it infinitesimal holonomy Lie algebra}.
The Lie subgroup of $\mathrm{GL}(T_pM,\R)$ generated by $\mathfrak{h}_p^\na$ is the {\it
infinitesimal holonomy group}
at $p$. The main result  is that if the infinitesimal holonomy group
has the
same dimension at every point $p$ of $M$ (which is the case when   $M$ and $\na$ are
analytic), then the restricted
holonomy group is equal to the infinitesimal holonomy group at every
point (see \cite{nomi}).
The linear connection $\na$ will be called invariant under parallelism in case $T^\na$
and $K^\na$ are both parallel with respect to $\na$.   The existence of  a linear
connection  $\na$ invariant under parallelism characterize (at lest locally) reductive
homogeneous spaces (see \cite{kostant}).  If $\na$ is invariant under parallelism then
\begin{equation}
 \label{eq01}\mathfrak{h}_p^\na=\left\{\; \sum K^\na(u_i,v_i),\; u_i,v_i\in T_pM\right\}.
\end{equation}

A vector field $A$ is an infinitesimal $\na$-transformation if and only if
for any couple of vector fields $X,Y$,
\begin{equation}\label{eqaffine}
  [A,\na_XY]=\na_{[A,X]}Y+\na_X[A,Y].
\end{equation}
On can see easily that this relation is equivalent to
\begin{equation}
 \label{eq2} \na_{X,Y}^2A+[\na_X,T_A^\na]Y=K^\na(X,A)Y,
\end{equation}where $\na_{X,Y}^2A=\na_X\na_YA-\na_{\na_XY}A$.

Let $\overline{\na}$ be another linear connections on  $M$. One
knows that $S=\overline{\na}-\na$ is a tensor field  of type $(1,2)$. By using a
terminology due to Kostant, we
will say that $\overline{\na}$ is rigid with respect to $\na$ whenever $S$ is
parallel with respect to $\na$. In this case, we have the following formula (see
\cite{kostant} Lemma 2):
\begin{equation}\label{eq3}
K^{\overline{\na}}(X,Y)=K^\na(X,Y)+[S_X,S_Y]+S_{T^\na(X,Y)}.\end{equation}
Let $G$ be a connected Lie group, $\G=T_eG$ its Lie algebra.
For any $u\in\G$ we denote by $u^+$ (resp. $u^-$) the left invariant (resp. the
right invariant) vector field associated to $u$.\\
 It is obvious that $G$ is a
reductive homogeneous space and hence, according to a result of Kostant (See
\cite{kostant} Theorem 2), $G$ admits a
linear connection invariant under parallelism. In fact $G$ admits many such connections
and we will use in this paper a special one, namely, the linear connection $\na^0$ given
by
$$\na_{u^+}^0v^+=\frac12[u^+,v^+],$$for any $u,v\in\G$.
 This connection is torsion free, invariant under parallelism, bi-invariant, complete and
its
curvature and holonomy Lie algebra are given by
\begin{equation}\label{eq4}K^{\na^0}(u^+,v^+)w^+=-\frac14[[u^+,v^+],w^+],
\quad u,v,w\in\G.
\end{equation}
\begin{equation}\label{eq5}
 \mathfrak{h}_e^{\na^0}=\ad_{[\G,\G]}.
\end{equation}A linear connection $\na$ on
$G$ is called bi-invariant if $\na$ is invariant by left and right multiplication.
The following lemma gives different characterizations of bi-invariant linear
connections on $G$.
\begin{Le}\label{lemma1}
 Let $\na$ be a linear connection on $G$. Then the following assertion are equivalent:
\begin{enumerate}
 \item $\na$ is a bi-invariant linear connection.
\item For any couple of left invariant (resp. right invariant) vector field $X,Y$,
$\na_XY$ is left invariant (resp. right invariant).
\item For any couple of left invariant  vector field $X,Y$,
$\na_XY$ is left invariant and the product $\al:\G\times\G\too\G$ given by
$\al(u,v)=(\na_{u^+}v^+)(e)$ satisfies
\begin{equation}\label{eqlemma1}[u,\al(v,w)]=\al([u,v],w)+\al(v,[u,w]).\end{equation}
\item For any couple of right invariant  vector field $X,Y$,
$\na_XY$ is right invariant and the product $\be:\G\times\G\too\G$ given by
$\be(u,v)=(\na_{u^-}v^-)(e)$ satisfies
$$[u,\be(v,w)]=\be([u,v],w)+\be(v,[u,w]).$$
\item $\na$ is left invariant and rigid with respect to $\na^0$.
\item $\na$ is right invariant and rigid with respect to $\na^0$.
\end{enumerate}

\end{Le}

{\bf Proof.} Since $G$ is connected, $\na$ is bi-invariant if and only if, for any
$u\in\G$, $u^+$ and $u^-$ are infinitesimal $\na$-transformations, i.e., according to
\eqref{eqaffine}, for any couple of vector fields $X,Y$,
$$
 [u^+,\na_XY]=\na_{[u^+,X]}Y+\na_X[u^+,Y]\esp
[u^-,\na_XY]=\na_{[u^-,X]}Y+\na_X[u^-,Y].
$$
Since $G$ is a parallelizable by left invariant vector field and these vector fields
commute with right invariant vector fields,
these equations are equivalent to
$$[u^+,\na_{v^+}w^+]=\na_{[u^+,v^+]}w^++\na_{v^+}[u^+,w^+]\esp
[u^-,\na_{v^+}w^+]=0,\;\quad v,w\in\G.$$
The group $G$ is also prallelizable by right invariant vector field and hence these
equations are also equivalent to
$$[u^-,\na_{v^-}w^-]=\na_{[u^-,v^-]}w^-+\na_{v^-}[u^-,w^-]\esp
[u^+,\na_{v^-}w^-]=0,\;\quad v,w\in\G.$$
On the other hand,   $\na$ is left
invariant and rigid with respect to $\na^0$ is equivalent to
$$[u^-,\na_{v^+}w^+]=0\esp \na^0_{u^+}(\na-\na^0)(v^+,w^+)=0.$$
Or
\begin{eqnarray*}
 \na^0_{u^+}(\na-\na^0)(v^+,w^+)&=&\na^0_{u^+}\left(\na_{v^+}w^+-\na_{v^+}^0w^+\right)
-\left(\na_{\na^0_{u^+}v^+}w^+-\na_{\na^0_{u^+}v^+}^0w^+\right)\\&&
-\left(\na_{v^+}\na^0_{u^+}w^+-\na_{v^+}^0\na^0_{u^+}w^+\right)\\
&=&\frac12[u^+,\na_{v^+}w^+]-\frac14[u^+,[v^+,w^+]]-\frac12
\na_{[{u^+},v^+]}w^+\\&&-\frac14[w^+,[u^+,v^+]]
-\frac12\na_{v^+}[u^+,v^+]-\frac14[v^+,[w^+,u^+]]\\
&=&\frac12[u^+,\na_{v^+}w^+]-\frac12
\na_{[{u^+},v^+]}w^+-\frac12\na_{v^+}[u^+,v^+].
\end{eqnarray*}
A similar computation holds when one replaces left invariant vector field by right ones.
Now we can get the desired equivalences easily. \hfill$\square$

\begin{rem}For any $u,v\in\G$,  $u^+(e)=u^-(e)$ and $[v^+,u^-]=0$, so  we get
$$(\na_{u^+}v^+)(e)=(\na_{u^-}v^+)(e)=(\na_{v^+}u^-)(e)=(\na_{v^-}
u^-)(e).$$Thus
$$\al(u,v)=\be(v,u).$$

\end{rem}

Let $\na$ be torsion free bi-invariant linear connection on $G$. As above, we define
$S=\na-\na^0$. It is clear that $S$ is bi-invariant and define a product $\circ$ on $\G$.
We have
$$u\circ v=\al(u,v)-\frac12[u,v]=\frac12\al(u,v)+\frac12\al(v,u)=
\frac12\al(u,v)+\frac12\be(u,v).$$
This product is obviously commutative and, according to Lemma \ref{lemma1} 3 and 4,
satisfies \eqref{eq0}. Since $\na$ is rigid with respect to $\na^0$, \eqref{eq3} holds and
can be written for any $u,v\in\G$,
$$K^{\na}(u,v)=K^{\na^0}(u,v)+[S_u,S_v].$$
Thus $K^{\na}=K^{\na^0}$ if and only if $[S_u,S_v]=0$ for any $u,v\in\G$, which is equivalent to   $\circ$ is associative.  Hence $\circ$ defines a Poisson structure on $\G$ if and only if
$\na$ and $\na^0$ have the same curvature. So we get the desired interpretation.
\begin{theo}\label{main}
Let $G$ be a  connected Lie group and $\G$ its Lie algebra. Then the following
assertions hold:
\begin{enumerate}
 \item Let  $\na$ be a left invariant linear connection on $\G$ an let $\circ$ the product
on $\G$
given by
$$u\circ v=(\na_{u^+}v^+)(e)-\frac12[u,v].$$Then $(\G,[\;,\;],\circ)$  is a  Poisson
algebra if and only if
$\na$ is torsion free, bi-invariant and has the same curvature as $\na^0$.
\item Let $\circ$ be a product on $\G$  such that $(\G,[\;,\;],\circ)$ is a Poisson
algebra. Then the linear connection on $G$
given by
$$\na_{u^+}v^+=\frac12[u^+,v^+]+(u\circ v)^+$$is torsion free, bi-invariant and has the
same
curvature as $\na^0$.
\end{enumerate}
\end{theo}
We call {\it special}  a torsion free bi-invariant linear connection which
has the same curvature as $\na^0$.\\

In  Riemannian geometry there is a notion of semi-symmetric spaces which is a direct generalization of locally symmetric spaces, namely, Riemannian manifolds for which the  curvature tensor $K$ satisfies $K.K=0$, i.e.,
\begin{equation}\label{semi} \na_X\na_YK-\na_Y\na_XK-\na_{[X,Y]}K=0, \end{equation}
for any vector fields $X,Y$. Semi-Riemannian symmetric spaces were investigated first by E. Cartan \cite{cartan} and studied by many authors (\cite{lichne2, couty,  szabo} etc). More generally, we call a torsion free linear connection on a manifold semi-symmetric if its  curvature tensor satisfies \eqref{semi}.
\begin{pr}\label{mainpr}Any special connection is semi-symmetric.
\end{pr}

{\bf Proof.} Let $\na$ be a special connection on a Lie group $G$. According to Theorem \ref{main}, its curvature $K$ is given by
\[ K(X,Y)Z=-\frac14[[X,Y],Z], \]for any left invariant vector fields $X,Y,Z$. Now, it was shown in \cite{szabo} pp. 532 that the equation \eqref{semi} is equivalent to
\[ [K(U,V),K(X,Y)]=K(K(U,V)X,Y)+K(X,K(U,V)Y), \]for any left invariant vector field $X,Y,U,V$. By replacing $K$ in this relation by its expression above we get the desired result.\hfill$\square$\\

Let $\na$ be a left invariant linear connection on $G$. The  holonomy Lie
algebra is the smallest subalgebra $\h^\na$ of $\mathrm{End}(\G)$ which contains all
$K^\na(u,v)$ and satisfying $[\na_u,\h^\na]\subset \h^\na$ for any $u\in\G$ (see
\cite{nomi}). It is clear
that it is difficult to compute $\h^\na$ explicitly. However,
the  holonomy algebra of a special linear connection  can be computed easily.
\begin{Le}
 \label{lemma2}Let $\na$ be a special connection on $G$. Then the holonomy Lie
algebra  of $\na$ is given by
$$\mathfrak{h}^\na_e=\ad_{[\G,\G]}+\mathrm{L}_{[[\G,\G],\G]}
=\ad_{[\G,\G]}+\mathrm{R}_{[[\G,\G],\G]},$$
where $\mathrm{L},\mathrm{R}:\G\too\mathrm{End}(\G)$ are given by $\mathrm{L}_u=\al(u,.)$
and $\mathrm{R}_u=\al(.,u)$ and $\al(u,v)=(\na_{u^+}v^+)(e)$.
\end{Le}
{\bf Proof.} This is a consequence of the following formulas which hold for any special
connection. We have, for any $u,v\in\G$,
\begin{eqnarray*}\;[\ad_u,\mathrm{L}_u]&=&\mathrm{L}_{[u,v]},\;
[\ad_u,\mathrm{R}_u]=\mathrm{R}_{[u,v]},\;[\mathrm{L}_u,\mathrm{L}_v]=\mathrm{L}_{[u,v]}
-\frac14\ad_{[u,v]},\\\; [\mathrm{R}_u,\mathrm{R}_v]&=&-\mathrm{R}_{[u,v]}
-\frac14\ad_{[u,v]}.\end{eqnarray*}These formulas will be stated rigorously in the next
section.
\hfill$\square$\\

A special connection which has also the same  holonomy Lie algebra as $\na^0$
is called {\it strongly special}.

\section{Poisson algebras and Poisson admissible algebras}\label{sectional}
In this section, we study Poisson algebras and Poisson admissible algebras in algebraic
point of
view, having in mind the results of the previous section.   \\

Let $(\G,[\;,\;])$ be a finite dimensional Lie algebra and $\al:\G\times\G\too\G$  a
product on $\G$. For any $u\in\G$,
we define $\mathrm{L}_u,\mathrm{R}_u:\G\too\G$ by
$$\mathrm{L}_uv=\al(u,v)\esp \mathrm{R}_uv=\al(v,u).$$
Suppose that $\al$ is Lie-admissible, i.e., for any $u,v\in\G$,
  $$\al(u,v)-\al(v,u)=[u,v].$$
Suppose also that $\al$ is bi-invariant, i.e., it satisfies \eqref{eqlemma1}.
It is obvious that the product $\circ$ on $\G$ given by
\begin{equation}\label{eq9}
u\circ v={\al}(u,v)-\frac12[u,v]=\frac12{\al}(u,
v)+\frac12 { \al } (v , u)
\end{equation}is bi-invariant which is equivalent to
$$[S_u,\ad_u]=S_{[u,v]},\eqno(*)$$for any $u,v\in\G$ where $S_uv=u\circ v$. If we denote
by $K^{{\al}}$ the curvature of $\al$, we get that
\begin{eqnarray*}
 K^{{\al}}(u,v)&:=&[\mathrm{L}_u,\mathrm{L}_v]-\mathrm{L}_{[u,v]}\\
&=&[S_u+\frac12\ad_u,S_v+\frac12\ad_v]-S_{[u,v]}-\frac12\ad_{[u,v]}\\
&\stackrel{(*)}=&[S_u,S_v]-\frac14\ad_{[u,v]}.
\end{eqnarray*}
This formula is the infinitesimal analog of \eqref{eq3}.
 Thus we have proved the following result.
\begin{pr}\label{prsal1}Let $\al:\G\times\G\too\G$ a Lie-admissible product on $\G$ and
 $\circ$ given by \eqref{eq9}. Then $(\G,[\;,\;],\circ)$ is a Poisson algebra
if and only if $\al$ is bi-invariant and, for any $u,v\in\G$,
 $$K^{{\al}}(u,v)=-\frac14\ad_{[u,v]}.$$
\end{pr}
In view of this proposition and Lemma \ref{lemma2}, we can introduce this definition.
\begin{Def}Let $\G$ be a finite dimensional Lie algebra.\begin{enumerate}
\item A  product $\al$ on  $\G$ is called {\it quasi-canonical} if it is
Lie-admissible,
bi-invariant and
 $$K^{{\al}}(u,v)=-\frac14\ad_{[u,v]}.$$
\item The holonomy Lie algebra of a quasi-canonical product $\al$ on $\G$ is the
subalgebra of the Lie
algebra
$\mathrm{End}(\G)$ given by
$$\mathfrak{h}^\al=\ad_{[\G,\G]}+\mathrm{L}_{[\G,[\G,\G]]}=
\ad_{[\G,\G]}+\mathrm{R}_{[\G,[\G,\G]]}.$$
\item A quasi-canonical  product $\al$ is called {\it strongly quasi-canonical} if
$\mathfrak{h}^\al=\ad_{[\G,\G]}$.\end{enumerate}
                                 \end{Def}

According to Proposition \ref{prsal1} there is a correspondence between the set of
Poisson products on a Lie algebra $\G$ and the set of quasi-canonical products on $\G$. We
call {\it strong} a Poisson product whose corresponding quasi-canonical product is
strongly quasi-canonical. The corresponding quasi-canonical product to the trivial Poisson
product
is $\al_0(u,v)=\frac12[u,v]$. \\

 Let $\al$ be a quasi-canonical product on a Lie algebra $\G$. Then \eqref{eqlemma1} is
equivalent to
\begin{equation} \label{eq6}
 [\ad_u,\mathrm{L}_v]=\mathrm{L}_{[u,v]},
\end{equation}or
\begin{equation} \label{eq7}
 [\ad_u,\mathrm{R}_v]=\mathrm{R}_{[u,v]},
\end{equation}for any $u,v\in\G$. Since  $\ad_u=
\mathrm{L}_u-\mathrm{R}_u,$
we get when we replace $\ad_u$ in \eqref{eq6} that
\begin{equation}\label{eq8}
 K^\al(u,v)=[\mathrm{R}_u,\mathrm{L}_v].
\end{equation}Note that the curvature of $\al$ vanishes if and only if $\al$ is
associative. In this case, the Lie algebra $\G$ is 2-nilpotent  because $K^\al(u,v)=-\frac14\ad_{[u,v]}$.\\

Let us define the infinitesimal analog of the notion of a connection which has
parallel curvature. A product $\al$ on a Lie algebra $\G$ has a parallel curvature if,
for any $u,v,w\in\G$,
$$\na K^\al(u,v,w):=[\mathrm{L}_u,K^\al(v,w)]-K^{\al}(\al(u,v),w)-K^\al(v,\al(u,w))=0.$$

Let us give now some properties of Poisson admissible algebras. The definition of a
Poisson admissible algebra was given in Section \ref{section1}.

\begin{pr}\label{prsal2} Let $(\mathrm{A},.)$ be an algebra. Then the following conditions
are
equivalent:
\begin{enumerate}
 \item $(\mathrm{A},.)$ is a Poisson admissible algebra.
\item For any
$u,v\in \mathrm{A}$,
$$[\mathrm{L}_u,\mathrm{L}_v]+[\mathrm{R}_u,\mathrm{R}_v]+2[\mathrm{L}_u,\mathrm{R}_v]=0
\esp K(u,v)=[\mathrm{R}_u , \mathrm { L }_v],$$where
$K(u,v):=[\mathrm{L}_u,\mathrm{L}_v]-\mathrm{L}_{[u,v]}.$
\item For any
$u,v\in \mathrm{A}$,
$$ [\mathrm{R}_u,\mathrm{R}_v]+\mathrm{L}_{[u,v]}+3[\mathrm{L}_u,\mathrm{R}_v]=0.
$$
\item For any
$u,v\in \mathrm{A}$,
$$ [\mathrm{L}_u,\mathrm{L}_v]-\mathrm{R}_{[u,v]}+3[\mathrm{R}_u,\mathrm{L}_v]=0.
$$

\end{enumerate}
\end{pr}

{\bf Proof.} For any $u,v\in \mathrm{A}$, put $\mathrm{L}_uv=\mathrm{R}_vu=u.v$,
$[u,v]=u.v-v.u$, $u\circ v=\frac12(u.v+v.u)$ and
$K(u,v)=[\mathrm{L}_u,\mathrm{L}_v]-\mathrm{L}_{[u,v]}$.\\

The algebra $(\mathrm{A},.)$ is  Poisson admissible   if and only if
$(\mathrm{A},[\;,\;])$ is a Lie algebra and $"."$ is quasi-canonical with respect to
$(\mathrm{A},[\;,\;])$. This is equivalent to
\begin{itemize}
\item $K(u,v)w+K(v,w)u+K(w,u)v=0,$  (Bianchi identity)
\item
$K(u,v)\stackrel{\eqref{eq8}}=[\mathrm{R}
_u , \mathrm { L } _v ] $ ,
\item $[\mathrm{L}_u+\mathrm{R}_u,\mathrm{L}_v+\mathrm{R}_v]=0.$ (The associativity of
$\circ$),
\end{itemize}for any $u,v\in \mathrm{A}$.
Since from the second condition we deduce that
$[\mathrm{R}_u , \mathrm { L }_v ]=-[\mathrm{R}_v , \mathrm { L }_u ]$ (the flexibility),
the conditions above are equivalent to
\begin{itemize}
\item $K(u,v)w+K(v,w)u+K(w,u)v=0,$
\item
$[\mathrm{L}_u,\mathrm{L}_v]-\mathrm{L}_{[u,v]}=[\mathrm{R}_u , \mathrm { L }_v ] $ ,
\item
$[\mathrm{L}_u,\mathrm{L}_v]+[\mathrm{R}_u,\mathrm{R}_v]+2[\mathrm{L}_u,\mathrm{R}_v]=0$.
\end{itemize}
So 1. implies 2. Now its obvious that 2. implies 3. \\
Let us show now that the third  condition implies the first one. Note first that
$$K(u,v)w=\mathrm{ass}(v,u,w)-\mathrm{ass}(u,v,w)\eqno(*)$$where
$$\mathrm{ass}(u,v,w)=(u.v).w-u.(v.w)=[\mathrm{R}_w,\mathrm{L}_u](v).\eqno(**)$$
Moreover, the third  condition implies that
$\mathrm{ass}(u,v,w)=-\mathrm{ass}(w,v,u)$ and hence the product $"."$ is Lie-admissible
if and only if
$$\mathrm{ass}(u,v,w)+\mathrm{ass}(v,w,u)+\mathrm{ass}(w,u,v)=0.\eqno(***)$$
Now from 3. we deduce that
$$3\mathrm{ass}(u,v,w)=
\mathrm{R}_u\circ\mathrm{R}_w(v)-\mathrm{R}_w\circ\mathrm{R}_u(v)+
\mathrm{R}_v\circ\mathrm{R}_w(u)-\mathrm{R}_v\circ\mathrm{R}_u(w)$$
and we get easily that "." is Poisson-admissible and consequently $K(u,v)=
[\mathrm{R}_u , \mathrm { L }_v]$. Now this relation and 3. implies that
$$[\mathrm{L}_u,\mathrm{L}_v]+[\mathrm{R}_u,\mathrm{R}_v]+2[\mathrm{L}_u,\mathrm{R}_v]
=0.$$ The condition 4.  is equivalent to 3. is a consequence of the following remark.
If $(\mathrm{A},.)$ is an algebra and $\star$ is the product given by
$u\star v=-v.u$ then $(\mathrm{A},.)$ is  Poisson admissible  if and only if
$(\mathrm{A},\star)$ is  Poisson  admissible. In this case the structures of Lie
algebras of $(\mathrm{A},.)$ and $(\mathrm{A},\star)$ coincident.
\hfill
$\square$\\

 Poisson admissible algebras is a subclass of Lie-admissible flexible algebras
studied in
\cite{BenOsb}. Recall that an algebra is called {\it flexible} if its associator satisfies
$$\mathrm{ass}(u,v,w)+\mathrm{ass}(w,v,u)=0,$$for any $u,v,w$. The third
characterization of Poisson admissible algebras in Proposition \ref{prsal2} appears first
in
\cite{MarRem}.

\begin{co}\label{co1} An associative algebra $(\mathrm{A},.)$ is Poisson admissible  if
and only if $(\mathrm{A},[\;,\;])$ is a 2-nilpotent Lie algebra.

\end{co}

An algebra $(\mathrm{A},.)$ is called $\mathrm{LR}$-algebra if, for any $u,v\in
\mathrm{A}$,
$$[\mathrm{L}_u,\mathrm{L}_v]=[\mathrm{R}_u,\mathrm{R}_v]=0.$$
It follows from Proposition \ref{prsal2} that a $\mathrm{LR}$-algebra is
Poisson admissible  if and only if it is associative.\\

Let us introduce now an important class of strong Poisson admissible algebras.
A {\it left Leibniz algebra} is an algebra $(\mathrm{A},.)$ such that for any
$u\in\mathrm{A}$, the left multiplication $\mathrm{L}_u$ is a derivation,
i.e., for any $v,w\in\mathrm{A}$,
$$ u.(v.w)=(u.v).w+v.(u.w).$$
This is equivalent to one of the two following relations
\begin{equation}\label{eqleibniz1}
 [\mathrm{L}_u,\mathrm{L}_v]=\mathrm{L}_{uv}\quad\mbox{or}\quad
[\mathrm{L}_u,\mathrm{R}_v]=\mathrm{R}_{uv}.
\end{equation}

A {\it right Leibniz algebra} is an algebra $(\mathrm{A},.)$ such that, for any
$u\in\mathrm{A}$,  the right
multiplication $\mathrm{R}_u$ is  a derivation,
i.e., for any $v,w\in\mathrm{A}$,
$$
(v.w).u=(v.u).w+v.(w.u).$$
This is equivalent to one of the two following relations
\begin{equation}\label{eqleibnizr}
 [\mathrm{R}_u,\mathrm{R}_v]=\mathrm{R}_{vu}\quad\mbox{or}\quad
[\mathrm{R}_u,\mathrm{L}_v]=\mathrm{L}_{vu}.
\end{equation}

An algebra which is left and right Leibniz is called {\it symmetric Leibniz algebra}.
Leibniz
algebras were introduced by Loday \cite{loday}. Many examples of symmetric Leibniz
algebras can be found in \cite{hidri}. By using \eqref{eqleibniz1} and
\eqref{eqleibnizr}, we get the following proposition.
\begin{pr}\label{prleibnizlr}
 The following assertions are equivalent:
\begin{enumerate}
 \item $(\mathrm{A},.)$ is a symmetric Leibniz algebra.
\item For any $u,v\in\mathrm{A}$,
$[\mathrm{L}_u,\mathrm{L}_v]=\mathrm{L}_{uv}=-\mathrm{R}_{uv}.$
\item For any $u,v\in\mathrm{A}$,
$[\mathrm{R}_u,\mathrm{R}_v]=\mathrm{R}_{vu}=-\mathrm{L}_{vu}.$

\end{enumerate}

\end{pr}

Any Lie algebra is a symmetric Leibniz algebra and any Leibniz algebra is Lie-admissible.
However, the class of symmetric Leibniz algebras contains strictly the class of Lie
algebras.
We can state now one of our main result.

\begin{theo}\label{theoleibniz}
 Let $(\mathrm{A},.)$ be a symmetric Leibniz algebra. Then $(\mathrm{A},.)$ is  Poisson
admissible, $"."$ is
strongly quasi-canonical on $(\mathrm{A},[\;,\;])$ and has parallel curvature.
\end{theo}
{\bf Proof.} Put
$$Q=[\mathrm{R}_u,\mathrm{R}_v]+\mathrm{L}_{[u,v]}+3[\mathrm{L}_u,\mathrm{R}_v].$$
By using  Proposition \ref{prleibnizlr}, we get
\begin{eqnarray*}
 Q&=&-\mathrm{L}_{vu}+\mathrm{L}_{uv}-\mathrm{L}_{vu}-3\mathrm{L}_{uv}
=0.
\end{eqnarray*}
Now, according to Proposition \ref{prsal2} we get that $(\mathrm{A},.)$ is  Poisson
admissible. On the other hand,
by using   Proposition \ref{prleibnizlr}, we
get
that, for any $u,v\in\mathrm{A}$,
$$\ad_{[u,v]}=2\mathrm{L}_{[u,v]}=4\mathrm{L}_{uv}\esp K(u,v)=\mathrm{L}_{vu}$$and hence
$\ad_{[\mathrm{A},\mathrm{A}]}=
\mathrm{L}_{[\mathrm{A},\mathrm{A}]}$. So the holonomy Lie algebra of
$"."$ is $\ad_{[\mathrm{A},\mathrm{A}]}$ which  prove that
$"."$
strongly quasi-canonical. Finally,
\begin{eqnarray*}
 \na K(u,v,w)&=&[\mathrm{L}_u,\mathrm{L}_{wv}]
-\mathrm{L}_{(uw)v}-\mathrm{L}_{w(uv)}=\\
&=&\mathrm{L}_{u(wv)}-\mathrm{L}_{(uw)v}-\mathrm{L}_{w(uv)}=0,
\end{eqnarray*}
which completes the proof.\hfill$\square$\\

By using the geometric interpretation  of Poisson structures introduced in Section
\ref{section2}, we get the following
interesting corollary.
\begin{co} Let $(\mathrm{A},.)$ be a real symmetric Leibniz algebra which is not a Lie
algebra and $G$ any connected Lie group associated to $(\mathrm{A},[\;,\;])$. Then the
left invariant connection on $G$ given by
$$\na_{u^+}v^+=(u.v)^+$$ is different from $\na^0$,  strongly special and its curvature
is parallel.

\end{co}

\begin{exem} We give here an example of a 4-dimensional real symmetric Leibniz algebra for
which we give the connected and simply connected Lie group associated to the underlying
Lie algebra and we give explicitly the two connections $\na^0$ and $\na$ appearing in the
corollary above.\\
We consider the symmetric Leibniz algebra product on $\R^4$ given in the canonical
basis $(e_1,e_2,e_3,e_4)$ by
$$e_1.e_1=e_4,\; e_2.e_1=e_3,\; e_3.e_1=e_4,\; e_1.e_2=-e_3,\; e_1.e_3=-e_4.$$
All the others products vanish.  One can check easily by using Proposition
\ref{prleibnizlr} that this product defines actually a symmetric Leibniz algebra.
The underlying Lie algebra say $\G=\R^4$ has its non-vanishing Lie brackets  given by
$$[e_1,e_2]=-2e_3\esp [e_1,e_3]=-2e_4.$$
It is a 3-nilpotent Lie algebra. The associated connected and simply connected Lie group
is $G=\R^4$ with the multiplication
given by Campbell-Baker-Hausdorff formula
$$xy=x+y+\frac12[x,y]+\frac1{12}[x,[x,y]]+\frac1{12}[y,[y,x]].$$
This formula can be written
\begin{eqnarray}xy&=&(x_1+y_1,x_2+y_2,x_3+y_3-(x_1y_2-x_2y_1),
x_4+y_4-(x_1y_3-x_3y_1)\nonumber\\&&+
\frac13x_1(x_1y_2-x_2y_1)+\frac13y_1(y_1x_2-y_2x_1))\label{example}.\end{eqnarray}
Recall that for any vector $u\in\G$, $u^+$ denotes the left invariant vector on $G$
associated to $u$. A straightforward computation using  \eqref{example} gives
\begin{eqnarray*}
 e_1^+&=&\frac{\partial}{\partial x_1}+x_2\frac{\partial}{\partial x_3}+(
x_3-\frac13x_1x_2)\frac{\partial}{\partial x_4},\\
e_2^+&=&\frac{\partial}{\partial x_2}-x_1\frac{\partial}{\partial
x_3}+\frac13x_1^2\frac{\partial}{\partial x_4},\\
e_3^+&=&\frac{\partial}{\partial x_3}-x_1
\frac{\partial}{\partial x_4},\; e_4^+=\frac{\partial}{\partial x_4},
\end{eqnarray*}where $(x_1,x_2,x_3,x_4)$ are the canonical coordinates of $\R^4$.
We consider the two torsion free linear connections on $G$ defined by the
formulas
\begin{equation}\na^0_{x^+}y^+=\frac12[x^+,y^+]\esp
\na_{x^+}y^+=(x.y)^+.\label{example1}
\end{equation}The dot here is the symmetric Leibniz product.
According to what above these two connections are bi-invariant, have the same curvature
and the same holonomy Lie algebra. Moreover, they  both have parallel curvature. Let us
compute the Christoffel symbols of $\na^0$ and $\na$ in the canonical coordinates
$(x_1,x_2,x_3,x_4)$. The computation is straightforward consisting of replacing $x$ and
$y$ in \eqref{example1} by $e_i,e_j$, $i,j=4,\ldots,1$. We get that the only non-vanishing
 Christoffel symbols are given by
$$\na^0_{\frac{\partial}{\partial x_1}}\frac{\partial}{\partial x_1}=
-\frac23x_2\frac{\partial}{\partial
x_4}\esp
\na^0_{\frac{\partial}{\partial x_1}}\frac{\partial}{\partial
x_2}=\frac13x_1\frac{\partial}{\partial x_4},$$
and
$$\na_{\frac{\partial}{\partial x_1}}\frac{\partial}{\partial x_1}=
(1-\frac23x_2)\frac{\partial}{\partial
x_4}\esp
\na_{\frac{\partial}{\partial x_1}}\frac{\partial}{\partial
x_2}=\frac13x_1\frac{\partial}{\partial x_4}.$$
We can also compute the exponential maps associated to $\na^0$ and $\na$ and we get
that $\exp_0:\G\too G$ is the identity, however $\exp:\G\too G$ is given by
$$\exp(a,b,c,d)=(a,b,c,d-\frac12a^2).$$

\end{exem}

\begin{pr}
 A left (right) Leibniz algebra  is  Poisson admissible if and only if it is  a symmetric
Leibniz algebra.
\end{pr}

{\bf Proof.} Suppose that $(\mathrm{A},.)$ is a left Leibniz algebra. According to
Proposition \ref{prsal2}, $(\mathrm{A},.)$ is  Poisson admissible if and only if,
for any
$u,v\in \mathrm{A}$,
$$ [\mathrm{L}_u,\mathrm{L}_v]-\mathrm{R}_{[u,v]}+3[\mathrm{R}_u,\mathrm{L}_v]=0.$$
From this relation and \eqref{eqleibniz1}, we get that $\mathrm{L}_{uv}=2
\mathrm{R}_{vu}+\mathrm{R}_{uv}$. On the other hand, \eqref{eqleibniz1} implies that
$\mathrm{L}_{uv}=-\mathrm{L}_{vu}$ so we deduce that $\mathrm{L}_{uv}=-\mathrm{R}_{uv}$
and we can achieve the proof by using Proposition \ref{prleibnizlr} and Theorem
\ref{theoleibniz}.\hfill $\square$\\

The following proposition gives a tool to build many symmetric  Leibniz algebras from old
ones.
\begin{pr}\label{prsal4}Let $\mathrm{A}$ be a symmetric Leibniz and $\mathrm{U}$ an
associative $\mathrm{LR}$-algebra then  $\mathrm{A}\otimes\mathrm{U}$
endowed with the product
$$(u\otimes a)(v\otimes b)=(uv)\otimes (ab)$$is a symmetric Leibniz algebra.

\end{pr}

{\bf Proof.} It is a straightforward computation.\hfill $\square$\\

We can state now our second main result.
\begin{theo}\label{theostrange}
Let $(\mathrm{A},.)$ be a Poisson admissible algebra and $\mathrm{U}$  an
associative $\mathrm{LR}$-algebra. Then the product on $\mathrm{A}\otimes\mathrm{U}$
given by
$$(u\otimes a)\star(v\otimes
b)=\frac12[u,v]\otimes(ab+ba)+\frac12u.v\otimes(3ab+ba)$$induces on
$\mathrm{A}\otimes\mathrm{U}$ a Poisson admissible algebra structure. Moreover, if $"."$
is
strongly quasi-canonical on $(\mathrm{A},[\;,\;])$ then $\star$ is also strongly
quasi-canonical
on $(\mathrm{A}\otimes\mathrm{U},[\;,\;])$.

\end{theo}

{\bf Proof.} Note first that since $\mathrm{U}$ is  an associative
$\mathrm{LR}$-algebra, for any $a_1,a_2,a_3\in\mathrm{U}$ and for any permutation
$\sigma$ of $\{1,2,3\}$,
$a_{\sigma(1)}a_{\sigma(2)}a_{\sigma(3)}=a_1a_2a_3$.\\
We will use Proposition  \ref{prsal2} and show that, for any
$u,v\in\mathrm{A}$ and $a,b\in\mathrm{U}$,
$$Q=[\mathrm{L}_{u\otimes a},\mathrm{L}_{v\otimes b}]-
\mathrm{R}_{[u\otimes a,v\otimes b]}+3[\mathrm{R}_{u\otimes a},\mathrm{L}_{v\otimes
b}]=0.$$
For any $w\in\mathrm{A}$ and $c\in\mathrm{U}$, we have
\begin{eqnarray*}
\;[\mathrm{L}_{u\otimes a},\mathrm{L}_{v\otimes b}](w\otimes c)&=&
 (u\otimes a)\star\left(\frac12[v,w]\otimes(bc+cb)+\frac12v.w\otimes(3bc+cb)\right)\\&&
-(v\otimes b)\star\left(\frac12[u,w]\otimes(ac+ca)+\frac12u.w\otimes(3ac+ca)\right)\\
&=&\left([u,[v,w]]+2u.[v,w]+2[u,vw]+4u(vw)\right)\otimes(abc)\\&&
-\left([v,[u,w]]+2v.[u,w]+2[v,uw]+4v(uw)\right)\otimes(abc)\\
&=&([u,[v,w]]+[v,[w,u]]+4[u,v].w+4[\mathrm{L}_{u},\mathrm{L}_{v}](w))\otimes(abc).
\end{eqnarray*}
Thus
$$[\mathrm{L}_{u\otimes a},\mathrm{L}_{v\otimes b}](w\otimes c)=
([u,[v,w]]+[v,[w,u]]+4[u,v].w+4[\mathrm{L}_{u},\mathrm{L}_{v}](w))\otimes(abc).$$
A similar computation gives
$$\mathrm{R}_{[u\otimes a,v\otimes b]}(w\otimes c)=
(4[w,[u,v]]+8w.[u,v])\otimes(abc),$$and
$$[\mathrm{R}_{u\otimes a},\mathrm{L}_{v\otimes
b}](w\otimes
c)=([[v,w],u]-[v,[w,u]]+2[[v,u],w]+4[\mathrm{R}_{u},\mathrm{L}_{v}](w))\otimes(abc).$$
By using Jacobi identity and the relation
$$[\mathrm{L}_u,\mathrm{L}_v]-\mathrm{R}_{[u,v]}+3[\mathrm{R}_u,\mathrm{L}_v]=0,$$
we get that $Q=0$ and hence $(\mathrm{A}\otimes\mathrm{U},.)$ is a Poisson admissible
algebra.\\
On the other hand, a direct computation gives, for any $u,v,w\in\mathrm{A}$ and any
$a,b,c\in\mathrm{U}$,
$$[[u\otimes a,v\otimes b],w\otimes c]=16[[u,v],w]\otimes(abc).$$
This shows that
$$[\mathrm{A}\otimes\mathrm{U},[\mathrm{A}\otimes\mathrm{U},\mathrm{A}\otimes\mathrm{U}]]
=[\mathrm{A},[\mathrm{A},\mathrm{A}]]\otimes \mathrm{U}^3,$$
and
$$\ad_{[u\otimes a,v\otimes b]}=16\ad_{[u,v]}\otimes\mathrm{L}_{ab}.$$
Moreover, one can check easily that for $u\in[\mathrm{A},[\mathrm{A},\mathrm{A}]]$ and
$a\in \mathrm{U}^3$,
$$\mathrm{L}_{u\otimes a}=(\ad_u+2\mathrm{L}_u)\otimes\mathrm{L}_a.$$
With all this formulas, one can show easily that if $"."$ is
strongly quasi-canonical on $(\mathrm{A},[\;,\;])$ then $\star$ is also strongly
quasi-canonical
on $(\mathrm{A}\otimes\mathrm{U},[\;,\;])$.
\hfill$\square$

\begin{pr}\label{prl}Let $(\G,[\;,\;])$ be a Lie algebra and $"."$ is a strongly
quasi-canonical product on $\G$. Then $\G^3=[\G,[\G,\G]]$ is two sided ideal of
$(\G,.)$,
$(\G^3,.)$ is a symmetric Leibniz algebra and the sequence
$$0\too(\G^3,.)\too (\G,.)\too (\G/\G^3,.)\too0$$is an exact sequence of Poisson
admissible algebras,
$(\G/\G^3,.)$ is associative and $(\G/\G^3,[\;,\;])$ is 2-nilpotent.
\end{pr}

{\bf Proof.} Since $"."$ is  strongly quasi-canonical then its holonomy Lie algebra is
equal to $\ad_{[\G,\G]}$ and hence $\mathrm{L}_{\G^3}\subset\ad_{[\G,\G]}$ and
$\mathrm{R}_{\G^3}\subset\ad_{[\G,\G]}$. Then for any $u\in\G$ and $v\in\G^3$ there
exists $w,t\in[\G,\G]$ such that $\mathrm{L}_u=\ad_w$ and $\mathrm{R}_u=\ad_t$. So
$u.v\in\G^3$, $v.u\in\G^3$ and $\mathrm{L}_u,\mathrm{R}_u$ are derivations of the
restriction of $"."$ to $\G^3$. We get that $\G^3$ is a two side ideal and $(\G^3,.)$ is
a symmetric Leibniz algebra. On the other hand, $(\G/\G^3,.)$ is a Poisson algebra and
$(\G/\G^3,[\;,\;])$ is 2-nilpotent so $(\G/\G^3,.)$ is associative.\hfill $\square$\\

In Proposition 23 of \cite{GozRem}, it was proved that there is no  non-trivial
Poisson
structure on a simple complex  Lie algebra. We finish this section by
generalizing this
result to any semi-simple Lie algebra over any field. We will show also
that in a perfect Lie algebra the canonical product
is the only strongly quasi-canonical product.

\begin{theo}\label{semisimple}
 \begin{enumerate}
  \item Let $\G$ be a perfect Lie algebra, i.e., $\G=[\G,\G]$. Then the  product
$\al_0(u,v)=\frac12[u,v]$ is the only strongly quasi-canonical product on $\G$.
\item Let $\G$ be a semi-simple Lie algebra. Then the  product $\al_0(u,v)=\frac12[u,v]$
is the only  quasi-canonical product on $\G$. In particular,   there is no non trivial
Poisson structure  on $\G$.
 \end{enumerate}

\end{theo}

{\bf Proof.} \begin{enumerate}
              \item Suppose that $"."$ is a strongly quasi-canonical product on $\G$ and
$[\G,\G]=\G$. We have shown in Proposition \ref{prl} that in this case the restriction of
$"."$ to $[\G,[\G,\G]]$
is a Leibniz product. Or $\G=[\G,[\G,\G]]$ and hence $(\G,.)$ is a Leibniz algebra. Now
from the relation $\mathrm{L}_{u.v}=-\mathrm{R}_{u.v}$ and the fact that $\G\G=\G$ we
deduce that $u.v=\frac12[u,v]$ for any $u,v\in\G$ and hence $"."$ is the canonical product
on $\G$.
\item Suppose that $"."$ is a  quasi-canonical product on a semi-simple Lie algebra $\G$,
denote by $\mathrm{L}_u$ and $\mathrm{R}_u$, respectively, the left and the right
multiplication by $u$ associated to $"."$ and put $S_u:=\mathrm{L}_u-\frac12\ad_u$. Note
first that since $\G$ is semi-simple, $\G=[\G,\G]$ and hence, by using \eqref{eq6}, we get
that for any $u\in\G$,
$$\tr(S_u)=0.$$
Consider
$$\mathcal{I}=\left\{u\in\G, S_u=0\right\}.$$
For any $u\in\G$ and any $v\in\mathcal{I}$, we have from \eqref{eq6} that
$$\mathrm{L}_{[u,v]}=[\ad_u,\mathrm{L}_v]=\frac12[\ad_u,\ad_v]=\frac12\ad_{[u,v]},$$and
hence $\mathcal{I}$ is an ideal of $(\G,[\;,\;])$. Let us show that
$\mathcal{I}=\G$. \\
Since $\G$ is semi-simple we have
$$\G=\oplus_{i=1}^{p}\G_i,$$ where $(\G_i)_{i=1}^{p}$ is a family of simple ideals of
$\G$, $$[\G_i,\G_i]=\G_i, \quad[\G_i,\G_j]=\{0\}\; \mbox{if}\; i\not=j.$$  We have for any
$i,j=1,\ldots,p$,
$$\G_i.\G_i\subset\G_i\esp \G_i.\G_j=\{0\}\quad\mbox{if}\; i\not=j.$$
Indeed, for any $u,v\in\G_i$ and for $j\not=i$ and $w\in\G_j$, we have
$$[w,u.v]=[w,u].v+u.[w,v]=0.$$ By using a similar argument, we get that if $i\not=j$,
$u\in\G_i$ and $v\in\G_j$, $u.v\in\G_i\oplus\G_j$. If $u=[a,b]$ with $a,b\in\G_i$, we get
$$u.v=[a,b].v=[a,b.v]\in\G_i.$$ Since $[\G_i,\G_i]=\G_i$ we get that $u.v\in\G_i$ and in
a similar way $u.v\in\G_j$ and hence $u.v=0$.\\
Suppose by contradiction that $\mathcal{I}\not=\G$. Since $\mathcal{I}$ is an ideal,
eventually by changing the indexation of the family $(\G_i)_{i=1}^{p}$, we can suppose
that there exists $1\leq r\leq p$ such that $$\G=
\mathcal{I}\oplus \mathcal{J}\esp \mathcal{J}=\oplus_{i=r}^p\G_i.$$
 For any $u\in\mathcal{J}$, we denote by
$\overline{S}_u$ the restriction of $S_u$ to $\mathcal{J}$. The product on $\mathcal{J}$
given by $u\circ v=
\overline{S}_uv$ is a Poisson product and hence it is commutative and associative. So, for
any
$u\in\mathcal{J}$, and any $n\in\N^*$,
$$\tr((\overline{S}_u)^n)=\tr(\overline{S}_{u^n})=\tr({S}_{u^n})=0,$$and hence
$\overline{S}_u$ is nilpotent. Since for any $u,v\in\mathcal{J}$,
$[\overline{S}_u,\overline{S}_v]=0$, we deduce by using Engel's Theorem that there exists
$u_0\in\mathcal{J}\setminus\{0\}$ such that $\overline{S}_u(u_0)=\overline{S}_{u_0}u=0$.
and hence $\overline{S}_{u_0}=0$. Since the restriction of $S_{u_0}$ to $\mathcal{I}$
vanishes, we deduce  that $S_{u_0}=0$ and hence $u_0\in
\mathcal{I}$ which constitutes a contradiction and achieves the proof.\hfill$\square$

             \end{enumerate}

\section{ Associative Poisson admissible    algebras }\label{ass}
We have shown in Corollary \ref{co1} that an associative algebra $(\mathrm{A},.)$ is
Poisson admissible if and only if
  $(\mathrm{A},[\;,\;])$ is
$2$-nilpotent, i.e., for any $u,v\in \mathrm{A}$,
\begin{equation}\label{eq12}\mathrm{L}_{[u,v]}=\mathrm{R}_{[u,v]}.\end{equation}
An associative algebra satisfying this condition will be called {\it associative Poisson
  admissible algebra}.
This  section is devoted to a description of such algebras.\\

Let $(\mathrm{A},.)$ be an associative Poisson admissible
algebra. We consider $$Z(\mathrm{A})=\left\{u\in\mathrm{A},
\mathrm{L}_u=\mathrm{R}_u\right\}.$$
Since $\mathrm{A}$ is associative, $Z(\mathrm{A})$ is a commutative associative subalgebra
of $\mathrm{A}$. Put
$$\mathrm{A}=\mathrm{V}\oplus Z(\mathrm{A}),$$where
$\mathrm{V}$ is a vectorial subspace of $\mathrm{A}$. According to this splitting, we get
that, for any $z\in Z(\mathrm{A})$ and $u,v\in\mathrm{V}$,
\begin{equation}
 z.u=u.z=\mathrm{P}_z(u)+\mathrm{Q}_u(z)\esp u.v=\mathfrak{a}(u,v)+\mathfrak{b}(u,v).
\end{equation}
These relations define two bilinear maps
$\mathfrak{a}:\mathrm{V}\times\mathrm{V}\too\mathrm{V}$,
$\mathfrak{b}:\mathrm{V}\times\mathrm{V}\too Z(\mathrm{A})$, and two linear maps
$\mathrm{Q}:\mathrm{V}\too\mathrm{End}( Z(\mathrm{A}))$,
$\mathrm{P}:Z(\mathrm{A})\too\mathrm{End}(\mathrm{V})$.
 The condition \eqref{eq12} is
equivalent to  $\mathfrak{a}$ symmetric and
$\mathfrak{b}(u,v)-\mathfrak{b}(v,u)=[u,v]$. The associativity is equivalent to the
following
relations:
\begin{enumerate}
 \item $\mathrm{P}_{zz'}=\mathrm{P}_{z}\circ\mathrm{P}_{z'}$,
$[\mathrm{Q}_u,\mathrm{L}_z](z')=\mathrm{Q}_{\mathrm{P}_{z'}(u)}(z)$,
\item  $\mathfrak{b}(\mathrm{P}_z(u),v)+\mathrm{Q}_{v}\circ{\mathrm{Q}_u(z)}=
\mathfrak{b}(u,\mathrm{P}_z(v))+\mathrm{Q}_{u}\circ{\mathrm{Q}_v(z)}=
\mathrm{Q}_{\mathfrak{a}(u,v))}(z)+z\mathfrak{b}(u,v),$
\item $\mathrm{P}_z(\mathfrak{a}(u,v))=\mathfrak{a}(u,\mathrm{P}_z(v))+
\mathrm{P}_{\mathrm{Q}_v(z)}(u),$
\item
$\mathfrak{b}(\mathfrak{a}(u,v),w)-\mathfrak{b}(u,\mathfrak{a}(v,w))=
\mathrm{Q}_{u}(\mathfrak{b}(v,w))-\mathrm{Q}_{w}(\mathfrak{b}(u,v)),$
\item $\mathfrak{a}(\mathfrak{a}(u,v),w)-\mathfrak{a}(u,\mathfrak{a}(v,w))
=\mathrm{P}_{\mathfrak{b}(v,w)}(u)-\mathrm{P}_{\mathfrak{b}(u,v)}(w).$
\end{enumerate}

So, we have shown that the associative and commutative algebra $Z(\mathrm{A})$, the vector
space $V$, and $\mathrm{P}$, $\mathrm{Q}$, $\mathfrak{a}$, $\mathfrak{b}$ as above
satisfying the conditions 1-5 describe entirely associative Poisson algebras.

\begin{pr}Let $(\G,[\;,\;])$ be a 2-nilpotent Lie algebra. Then there is on $\G$ a
quasi-canonical product different from the canonical one.

\end{pr}

{\bf Proof.} Put $\G=Z(\G)\oplus V$ and consider the product on $\G$ given, for any
$z,z'\in Z(\G),$ $u,v\in V$, by
$$z.z'=u.z=z.u=0 \esp u.v=\mathfrak{s}(u,v)+\frac12[u,v],$$
where $\mathfrak{s}$ is any non trivial symmetric bilinear map from $V\times V$ to
$Z(\G)$. It is easy to check that this product is quasi-canonical and different from the
canonical one.\hfill$\square$

\section{Symplectic Poisson algebras}\label{sectionsp}
In this section, we study an important class of Poisson algebras. To introduce theses
algebras we recall some classical results on symplectic Lie groups and introduce a new
symplectic linear connection.\\

Let $(G,\Om)$ be a symplectic Lie group, i.e., a Lie group $G$ endowed with a left
invariant symplectic form $\Om$. It is well-known  that the linear
connection given by
the formula
\begin{equation}\label{borde}\Om({\na}^{\mathrm{a}}_{u^+}v^+,w^+)=-\Om(v^+,[u^+,w^+]),
\end{equation}where $u,v,w\in\G$,  defines a left invariant flat and torsion free
connection ${\na}^{\mathrm{a}}$. Moreover, ${\na}^{\mathrm{a}}\Om$  never vanishes unless
$G$ is abelian. So we can define a tensor field $\mathrm{N}$ by the relation
$$\na^{\mathrm{a}}_{u^+}\Om(v^+,w^+)=\Om(\mathrm{N}(u^+,v^+),w^+).$$
The linear connection given by
$$\na^{\mathrm{s}}_{u^+}v^+=\na^{\mathrm{a}}_{u^+}v^++\frac13\mathrm{N}(u^+,v^+)+
\frac13\mathrm{N}(v^+,u^+)$$ is left invariant torsion free and symplectic, i.e.,
$\na^{\mathrm{s}}\Om=0$. This construction follows a general scheme which permit to build
symplectic connection from any connection (see \cite{cahen}). A straightforward
computation gives that $\na^{\mathrm{s}}$ can be defined by the following formula
\begin{equation}\label{eqsymplectic}
\Om(\na^{\mathrm{s}}_{u^+}v^+,w^+)=\frac13\Om([u^+,v^+],w^+)+
\frac13\Om([u^+,w^+],v^+).
\end{equation}
This formula shows that on any symplectic Lie group there exists a canonical torsion free
symplectic  connection.

Let $(\G,\om)$ be  the Lie algebra of $G$ endowed with the value of $\Om$ at $e$.
We denote by $\al^{\mathrm{a}}$ and $\al^{\mathrm{s}}$ the product on $\G$ induced,
respectively, by
$\na^{\mathrm{a}}$ and $\na^{\mathrm{s}}$. We have, for any $u,v\in\G$,
\begin{equation}\label{eqnewconnection}
 \al^{\mathrm{a}}(u,v)=-\ad_u^*v\esp \al^{\mathrm{s}}(u,v)=
\frac13\left(\ad_uv-\ad_u^*v\right),
\end{equation}where $\ad_u^*$ is the adjoint of $\ad_u$ with respect to $\om$.\\

Conversely, given a symplectic Lie algebra, the formulas
\eqref{eqnewconnection} defines on $\G$ two Lie-admissible products whose one is left
symmetric and the other one is symplectic. Let us see under which conditions these
products
are quasi-canonical.

\begin{pr}\label{prsp1} Let $(\G,\om)$ be a symplectic Lie algebra and
$\al^{\mathrm{a}}$, $\al^{\mathrm{s}}$ the product given by \eqref{eqnewconnection}. Then
the following
assertions are equivalent:
\begin{enumerate}
 \item $\al^{\mathrm{a}}$ is quasi-canonical.
\item $\al^{\mathrm{s}}$ is quasi-canonical.
\item $\G$ is 2-nilpotent and, for any $u,v\in\G$, $[\ad_u,\ad_v^*]=0$.
\end{enumerate}
Moreover, if one of the conditions above holds then $(\G,\al^{\mathrm{a}})$ and
$(\G,\al^{\mathrm{s}})$ are both
associative $\mathrm{LR}$-algebras.

\end{pr}

{\bf Proof.} Not first that $K^{\al^{\mathrm{a}}}=0$ and the  left and
right multiplications associated to $\al^a$ are given by
$$\mathrm{L}_u^{\mathrm{a}}=-\ad_u^*\esp \mathrm{R}_u^{\mathrm{a}}=-\ad_u^*-\ad_u.$$ The
product
$\al^{\mathrm{a}}$ is quasi-canonical if and only if, for any $u,v\in\G$,
$$K^{\al^{\mathrm{a}}}(u,v)=[\mathrm{R}_u^{\mathrm{a}},\mathrm{L}_u^{\mathrm{a}}]\esp
[\mathrm{L}_u^{\mathrm{a}}+\mathrm{R}_u^{\mathrm{a}},\mathrm{L}_v^{\mathrm{a}}+
\mathrm{R}_v^{\mathrm{a}}] =0$$which is obviously
equivalent to
$$\ad_{[u,v]}=[\ad_u,\ad_v^*]=0.$$
On the other hand, we have
$$\mathrm{L}_u^{\mathrm{s}}=\frac13(\ad_u-\ad_u^*)\esp
\mathrm{R}_u^{\mathrm{s}}=-\frac13(2\ad_u+\ad_u^*),$$and hence

\begin{eqnarray*}
 K^{\al^{\mathrm{s}}}(u,v)&=&[\mathrm{L}_u^{\mathrm{s}},\mathrm{L}_v^{\mathrm{s}}]-\mathrm
{L}_{[u,v]}^{\mathrm{s}}\\
&=&\frac19\left([\ad_u,\ad_v]-[\ad_u,\ad_v^*]-[\ad_u^*,\ad_v]+[\ad_u^*,\ad_v^*]\right)\\&&
-\frac13\ad_{[u,v]}+\frac13\ad_{[u,v]}^*,\\
\;[\mathrm{R}_u^{\mathrm{s}},\mathrm{L}_v^{\mathrm{s}}]&=&-\frac19
\left(2[\ad_u,\ad_v]-2[\ad_u,\ad_v^*]+[\ad_u^*,\ad_v]-[\ad_u^*,\ad_v^*]\right).
\end{eqnarray*}
Thus $K^{\al^{\mathrm{s}}}(u,v)=[\mathrm{R}_u^{\mathrm{s}},\mathrm{L}_v^{\mathrm{s}}]$ if
and only if
\begin{equation}
 [\ad_u,\ad_v^*]=\ad_{[u,v]}^*.
\end{equation}
Let us compute
$$Q=[\mathrm{L}_u^{\mathrm{s}}+\mathrm{R}_u^{\mathrm{s}},
\mathrm{L}_v^{\mathrm{s}}+\mathrm{R}_v^{\mathrm{s}}].$$
We have
\begin{eqnarray*}
 Q&=&\frac19\left([\ad_u,\ad_v]+2[\ad_u,\ad_v^*]+2[\ad_u^*,\ad_v]+4[\ad_u^*,\ad_v^*]
\right).
\end{eqnarray*}
So we get that $K^{\al^{\mathrm{s}}}(u,v)=[\mathrm{R}_u^{\mathrm{s}},
\mathrm{L}_v^{\mathrm{s}}]$  and
$Q=0$ if and only if
$$[\ad_u,\ad_v^*]=\ad_{[u,v]}=0.$$
\hfill$\square$

A {\it  symplectic Poisson algebra} is a 2-nilpotent symplectic Lie algebra $(\G,\om)$
satisfying, for any $u,v\in\G$,
\begin{equation}\label{eqsp}
[\ad_u,\ad_v^*]=0.
\end{equation}

\begin{pr} Let $(\G,\om)$ be 2-nilpotent symplectic Lie algebra which carries a
bi-invariant pseudo-Euclidean product $B$. Then $(\G,\om)$ is a symplectic Poisson
algebra.

\end{pr}

{\bf Proof.} We consider the isomorphism of $\G$ given by
$$\om(u,v)=B(Du,v).$$
It is easy to check by using the fact that $B$ is bi-invariant and $\om$ is symplectic
that $D$ is derivation of $\G$ and that, for any $u\in\G$,
$\ad_u^*=-D^{-1}\circ\ad_u\circ D.$ Now, for any $u,v\in\G$,
\begin{eqnarray*}
 [\ad_u^*,\ad_v]&=&\ad_v\circ D^{-1}\circ\ad_u\circ D-D^{-1}\circ\ad_u\circ D\circ\ad_v.
\end{eqnarray*}
Since $D\circ\ad_v=\ad_v\circ D+\ad_{Dv}$ and $\G$ is 2-nilpotent we get that
$$D^{-1}\circ\ad_u\circ D\circ\ad_v=0.$$ On the other hand, $D^{-1}[\G,\G]=[\G,\G]$ so we
get since $\G$ is 2-nilpotent that $\ad_v\circ D^{-1}\circ\ad_u\circ D=0$ and finally,
$[\ad_u^*,\ad_v]=0$ which show that $(\G,\om)$ is  a symplectic Poisson algebra.
\hfill$\square$

\begin{exem} Let  $(\G,[~,~])$  be a $2-$nilpotent Lie algebra. Then  $\G=
{\mathcal V}\oplus Z({\G}),$ where  ${\cal V}$ is a vector subspace  of  $\G$ such that
$[{\mathcal V},{\mathcal V}]\subseteq Z(\G).$
The endomorphism $D$ of  $\G$ defined by: $$D(v)= v\esp  D(z)= 2z,
\quad\mbox{for all}\; v \in {\mathcal V}, z\in Z({\G}),$$ is an invertible derivation of
$\G.$\\
Now , the vector space ${\mathcal G}= {\G}\oplus{\G}^*$ endowed with the following
product:
$$[u+\al,v+\be]= [u,v] +\al\circ\ad _{v}- \be\circ\ad _{u}  ,\quad \mbox{for all}\;
u,v \in {\G}, \al, \be \in {\G}^*,$$
is a  $2-$nilpotent Lie algebra. Moreover, the bilinear form  $B: {\mathcal G}\times
{\mathcal G} \rightarrow {\Bbb K}$ defined  by:
$$ B(u+\al,v+\be)= \al(v) + \be(u), \quad \mbox{for all}\; u,v \in {\G},  \al, \be \in
{\G}^*,
$$
is non-degenerate, bi-invariant and symmetric. Then $({\mathcal G},B)$ is a $2-$nilpotent
quadratic Lie algebra. An easy computation shows that the endomorphism  $\de$ of
$\mathcal G$ defined by:
$$\de(u)= D(u)~~\mbox{and}~~\de(\al)= -\al\circ D,\quad \mbox{for all}\; u\in {\G},
\al\in
{\G}^*,$$
is an invertible derivation of  $\mathcal G$ which is skew-symmetric with respect to $B$.
Consequently,  the bilinear form  defined by: $$\om(X,Y)= B(\de(X),Y),\quad \mbox{for
all}\;X,Y \in {\mathcal G},$$ is a symplectic
structure on  $\mathcal G$. Finally, $({\mathcal G},\om)$ is symplectic Poisson
algebra.\end{exem}
Let us give now the inductive description of symplectic Poisson algebras.
Let  $({\G},[~,~]_{\G},\om)$  be  a symplectic Poisson algebra. Since $\G$ is
nilpotent Lie algebra,
according to   \cite{DarMed},  $({\G},[~,~]_{\G},\om)$ is the symplectic double extension
of a symplectic Lie algebra
$({\h},[~,~]_{\h},\overline{\om})$ of dimension  $\dim{\G}-2$ by the one-dimensional Lie
algebra by means an element $(D,z)$ of $\mathrm{Der}({\mathfrak{h}})\times
{\mathfrak{h}}$.  This means that  ${\G}=
{\Bbb K}e \oplus {\h}\oplus {\Bbb K}d$ and
\begin{enumerate}
\item for any $a,b\in\h$,
\begin{eqnarray*}\;[a,b]_{\G}&=& [a,b]_{\h} + \overline{\om}((D+D^*)(a),b)e,\;
[d,d]_{\G}=0,\\\; [a,d]_{\G}&=& D(a)+\overline{\om}(z,a)e,\;
[e,{\G}]_{\G}=\{0\},\end{eqnarray*} where  $D^*$ the
adjoint of $D$ with respect to $\overline{\om},$
\item  $\om_{\vert_{{\h}\times {\h}}}=\overline{\om}, \om(e,d)=1, \om(e,{\h })=
\om(d,{\h})= \{0\}.$
\end{enumerate} The fact that $(\G,\om)$ is symplectic Poisson algebra is equivalent
to
$$\ad_{[u,v]}=[\ad_u,\ad_v^*]=0,$$for any $u,v\in\G$.

The first condition which means that $\G$ is 2-nilpotent Lie algebra is equivalent to:
\begin{itemize}
 \item $\h$ is a 2-nilpotent Lie algebra,
\item $D(\h)\subset Z(\h)$,
\item $D_{|[\h,\h]_{\h}}=D_{|[\h,\h]_{\h}}^*=0,$ $\overline{\om}([\h,\h]_{\h},z)=0$,
\item $D^2=D^*\circ D=0$ and $D^*(z)=0$.
\end{itemize}
Let us compute $\ad_u^*$ for any $u\in\G$. A straightforward computation gives, for any
$a,b\in\h$,
\begin{eqnarray*}\ad_{a}^*b&=&\ad_{a}^{\h*}b+\overline{\om}(b,D(a))e,\;\ad_{a}
^*d=-(D+D^*)(a)+\overline { \om}(a,z)e,\\
\ad_{a}^*e&=&\ad_{d}^*e=0,\;\ad_{d}^*d=-z,\;\ad_{d}^*a=-D^*(a).
\end{eqnarray*}So
\begin{eqnarray*}
 \;[\ad_a,\ad_b^*](c)&=&[a,\ad_{b}^{\h*}c]-\ad_{b}^{*}[a,c]=[a,\ad_{b}^{\h*}c]-\ad_{b}^{
\h* }[a,c]_{\h}+\overline{\om}([a,c]_{\h},D(b))e\\
&=&[\ad_{a}^{\h},\ad_{b}^{\h*}](c)+\overline{\om}((D+D^*)(a),\ad_{b}^{\h*}c)e
+\overline{\om}([a,c]_{\h},D(b))e,\\
\;[\ad_a,\ad_b^*](d)&=&-[a,(D+D^*)(b)]-\ad_b^*D(a)\\
&=&-[a,(D+D^*)(b)]_\h-\ad_{b}^{\h*}D(a)\\&&
-\overline{\om}((D+D^*)(a),(D+D^*)(b))e-\overline{\om}(D(a),D(b))e,\\
\;[\ad_a,\ad_d^*](b)&=&-[a,D^*(b)]+D^*([a,b]_\h)=
-[a,D^*(b)]_\h\\&&-\overline{\om}((D+D^*)(a),D^*(b))e+D^*([a,b]_\h),\\
\;[\ad_a,\ad_d^*](d)&=&-[a,z]+D^*\circ D(a)\\
&=&-[a,z]_\h-\overline{\om}((D+D^*)(a),z)e+D^*\circ D(a),\\
\;[\ad_d,\ad_d^*](a)&=&-[d,D^*(a)]-\ad_d^*[a,d]=D\circ D^*(a)
-\overline{\om}(z,D^*(a))e+D^*\circ D(a),\\
\;[\ad_d,\ad_d^*](d)&=&-D(z).\\
\end{eqnarray*}
From these relations, we get that $(\G,\om)$ is symplectic  Poisson algebra if and
only if
\begin{itemize}
 \item $(\h,\overline{\om})$ is a symplectic Poisson algebra,
\item $D(\h)\subset Z(\h)$, $D^*(\h)\subset Z(\h),$
\item $D^2=D^*\circ D=D\circ D^*=0$, $D^*(z)=D(z)=0$ and $z\in
Z(\h)$.
\end{itemize}

 An element $(D,z)$ of
 $\mathrm{Der}({\h})\times {\h}$ which verifies the conditions  above
will be called  {\it admissible}.

To summarize, we have proved the following theorem.

\begin{theo}\label{theods}
Let  $({\G},\om)$  be a  symplectic   Lie algebra.
Then $({\G},\om)$ is a symplectic Poisson algebra if and only if   it is a
symplectic double extension of a symplectic  Poisson algebra
$({\h},\overline{\om})$ of
dimension $\dim{\G}-2$ by the one dimensional Lie algebra by means of  an  admissible
element $(D,z)\in\mathrm{Der}({\h})\times {\h}$.
\end{theo}

There is only one 2-dimensional symplectic Poisson algebra, namely the two
2-dimensional abelian Lie algebra $\h_0$ endowed with a symplectic form $\om_0$. There
exists a basis $\B=\{e_1,e_2\}$ of $\h_0$ such that $\om_0(e_1,e_2)=1$. An element
$(D,z)\in\mathrm{Der}({\h_0})\times {\h_0}$ is admissible if and only if $z$ is any
element of ${\h_0}$ and the matrix of $D$ in the basis $\B$ has one of the following forms
$$\left(\begin{array}{cc}0&a\\0&0\end{array}\right),
\left(\begin{array}{cc}0&0\\a&0\end{array}\right),
\left(\begin{array}{cc}a&b\\-\frac{a^2}b&-a\end{array}\right),\;b\not=0.$$
So we get all four dimensional symplectic Poisson algebras.
\begin{pr}Let $\G$ be a 4-dimensional Lie algebra. Then $\G$ is a symplectic Poisson
algebra iff it is isomorphic to one of the following symplectic Lie algebras:
\begin{enumerate}
\item $\mathrm{span}\{e,e_1,e_2,d\}$ with the non vanishing brackets
$$[e_1,d]=-z_2e,\;[e_2,d]=-ae_1+z_1e,$$and the symplectic form satisfying
$$\om(e,d)=\om(e_1,e_2)=1,\;\om(e,e_1)=\om(e,e_2)=\om(d,e_1)=\om(d,e_2)=0.$$
 \item $\mathrm{span}\{e,e_1,e_2,d\}$ with the non vanishing brackets
$$[e_1,d]=ae_1-\frac{a^2}be_2-z_2e,\;[e_2,d]=be_1-ae_2+z_1e,$$and the symplectic form
satisfying
$$\om(e,d)=\om(e_1,e_2)=1,\;\om(e,e_1)=\om(e,e_2)=\om(d,e_1)=\om(d,e_2)=0.$$

\end{enumerate}

\end{pr}
We finish this section by giving an important geometric property of real symplectic
Poisson algebras. \\
Let $(\G,\om)$ be a non abelian real symplectic Poisson algebra and $G$ a connected Lie
group having $\G$ as its Lie algebra. The symplectic form $\om$ defines on $G$ a
symplectic left invariant form $\Om$. Consider the two linear connections
$\na^{\mathrm{a}}$ and
$\na^{\mathrm{s}}$  defined on $G$ by \eqref{borde}-\eqref{eqsymplectic}. These two
connections are bi-invariant, flat, complete and $\na^{\mathrm{s}}\Om=0$. It was shown in
\cite{boumedina} that $\Om$ is polynomial of degree at most $\dim G-1$ in any affine
coordinates chart associated to $\na^{\mathrm{a}}$. The following result gives a more
accurate statement on the polynomial nature of $\Om$.
\begin{theo}\label{polynomial}
With the hypothesis and the notations above we have
$$(\na^{\mathrm{a}})^3\Om=0.$$In particular, $\Om$ is
polynomial of degree at  least one and at most $2$ in any affine
coordinates chart associated to $\na^{\mathrm{a}}$. Moreover, if the restriction of $\om$
to $[\G,\G]$ does not vanish then the degree is 2.

\end{theo}

{\bf Proof.}  For any $u,v,x,y\in\G$, an easy computation gives
$$\na^{\mathrm{a}}_{u^+}\Om(x^+,y^+)=\Om(u^+,[x^+,y^+]),$$
and hence
$$\na^{\mathrm{a}}_{u^+}\na^{\mathrm{a}}_{v^+}\Om(x^+,y^+)
=\Om([\na^{\mathrm{a}}_{u^+}x^+,y^+]+[x^+,\na^{\mathrm{a}}_{u^+}y^+],v^+).$$
Now since $\na^{\mathrm{a}}$ is bi-invariant then
$$[\na^{\mathrm{a}}_{u^+}x^+,y^+]+[x^+,\na^{\mathrm{a}}_{u^+}y^+]=
\na^{\mathrm{a}}_{[u^+,y^+]}x^+
+\na^{\mathrm{a}}_{[x^+,u^+]}y^++2\na^{\mathrm{a}}_{u^+}[x^+,y^+].$$
By using \eqref{borde} and the fact that $\G$ is 2-nilpotent, we get
$$\na^{\mathrm{a}}_{u^+}\na^{\mathrm{a}}_{v^+}\Om(x^+,y^+)=
2\Om([x^+,y^+],[u^+,v^+]).$$
By using the same arguments as above one can get easily that $(\na^{\mathrm{a}})^3\Om=0.$
The properties of the degree of $\Om$ are an immediate consequence of formulas above.
\hfill$\square$\\

It was proved in \cite{gold2} that a compact affine manifold $M$ has a polynomial
Riemannian metric iff  $M$ is finitely covered by a complete affine nilmanifold. An
affine nilmanifold is  of the form $\Ga/N$ where $N$ is a simply-connected nilpotent Lie
group with a left invariant affine structure and $\Ga$ is a discrete subgroup of $N$.
According to the results of this section, if $G$ is the simply-connected Lie group
associated to a non abelian symplectic Poisson Lie algebra and $\Ga$ is a co-compact
discrete subgroup of $G$ then $\Ga/G$ is a compact nilmanifold which carries two affine
structures and a symplectic form which is parallel for one affine structure and polynomial
of degree at least 1 and at most 2 for the other one. It is natural to ask if there is a
symplectic analog of Goldman's Theorem in \cite{gold2}.

\section{Metrizability of special connections}\label{metri}

In this section we study the problem of metrizability of special connections on Lie
groups. Given a connected Lie group $G$ with  $\na$ a special connection, does
exist
on $G$ a left invariant  pseudo-Riemannian metric whose associated Levi-Civita connection
is $\na$ ? Remark that if such a metric exists and it is bi-invariant then $\na$ coincides
with $\na^0$. The following proposition gives an answer to this question when the metric
is Riemannian.

\begin{pr}
 Let $\G$ be a real  Lie algebra and $\prs$ an Euclidean product on $\G$ such that the
associated Levi-Civita product is quasi-canonical. Then $\prs$ is
bi-invariant and hence the Levi-Civita product coincides with the canonical product.
\end{pr}
{\bf Proof.} We have $$\G=[\G,\G]\oplus[\G,\G]^\perp.$$
 Since for any $u,v\in\G$, $K(u,v)=-\frac14\ad_{[u,v]}$ and $K(u,v)$ is
skew-symmetric, we deduce that, for any $w\in[\G,\G]$, $\ad_w$ is skew-symmetric.
From this remark and the relation
$$2\langle u.v,w\rangle=\langle[u,v],w\rangle+
\langle[w,v],u\rangle+\langle[w,u],v\rangle,$$one can deduce easily that,
for any $u,v\in[\G,\G]$ and any $x,y\in[\G,\G]^\perp$,
$$u.v=\frac12[u,v]\esp x.y=\frac12[x,y].$$
Let $u\in[\G,\G]$ and $v\in[\G,\G]^\perp$, since $\ad_u$ is skew-symmetric, we get for
any $w\in\G$,
$$\langle[u,v],w\rangle=-\langle v,[u,w]\rangle=0,$$and hence
$u.v=v.u$. Moreover, for any $x\in[\G,\G]$,
$$\langle u.v,x\rangle=-\langle v,u.x\rangle=\frac12\langle v,[x,u]\rangle=0.$$
Thus $u.v=v.u\in[\G,\G]^\perp$. Now
\begin{eqnarray*}
 2\langle u.v,u.v\rangle&=&\langle[u,v],u.v\rangle+\langle[u.v,u],v\rangle+
\langle[u.v,v],u\rangle\\
&=&0,
\end{eqnarray*}since $[u ,v]=0$ and $[v,u.v]=[v,u].v+u.[v,v]=0$. Thus $u.v=0$ which
completes the proof.
\hfill $\square$\\

The proposition above is not true in general when the $\prs$ is not positive definite. We give now a description of  all real  Lie algebras endowed with  an pseudo-Euclidean product  such that the
associated Levi-Civita product is quasi-canonical and the derived ideal is non-degenerate.\\

 Consider $(\h,\prs_0)$ a Lie algebra endowed with a
bi-invariant pseudo-Euclidean product. Let $(V,B)$ be a vector space with a nondegenerate
symmetric bilinear form. We can split $V=V_0\oplus U\oplus\overline{V}_0$ such that the
restriction of $B$ to $U$ is positive definite and the map $
V_0\times \overline{V}_0\too\R$, $(u,v)\too B(u,v)$ is non-degenerate. Finally, consider
any bilinear skew-symmetric map $\ga:\overline{V}_0\times\overline{V}_0\too Z(\h)$. We
consider now $\G=\h\oplus V$ endowed with $\prs=\prs_0+B$ and the bracket
for which $V_0\oplus U\subset Z(\G)$, $[\h,V]=0$, the restriction to $\h$ coincides with
the
initial bracket and for any $u,v\in \overline{V}_0$, $[u,v]=\ga(u,v)$. Then one can check
easily that the Levi-Civita product of $\prs$ is quasi-canonical and $\prs$ is not
bi-invariant. By a direct computation we can see easily that the  curvature tensor is parallel which gives examples of locally symmetric pseudo-Riemannian spaces.

Sa\"id Benayadi\\
Universit\'e de Lorraine,  IECL, CNRS-UMR 7502,\\ Ile du Saulcy, F-57045 Metz
cedex
1, France.\\
e-mail: said.benayadi@univ-lorraine.fr\\

\noindent Mohamed Boucetta\\
Universit\'e Cadi-Ayyad\\
Facult\'e des sciences et techniques\\
BP 549 Marrakech Maroc\\
e-mail: boucetta@fstg-marrakech.ac.ma

\end{document}